%(2) expari(al)=as
%(3) STARTING REAL bold in appendix C?  proof incomplete?
\magnification=\magstep1
\overfullrule=0pt
\input amssym
\def\ARI{{\rm{ARI}}}

\def\GARI{{\rm{GARI}}}
\def\GVARI{{\rm{\overline{GARI}}}}
\def\ari{{\rm{ari}}}

\def\ganit{{\rm{ganit}}}
\def\dur{{\rm{dur}}}
\def\dar{{\rm{dar}}}
\def\adari{{\rm{Ad_{\ari}}}}
\def\logari{{\rm{log_{ari}}}}
\def\gari{{\rm{gari}}}
\def\invgari{{\rm{inv}_{gari}}}
\def\mu{{\rm{mu}}}
\def\lu{{\rm{lu}}}
\def\neg{{\rm{neg}}}
\def\mantar{{\rm{mantar}}}
\def\expari{{\rm{exp_{ari}}}}

\def\preari{{\rm{preari}}}
\def\push{{\rm{push}}}
\def\amit{{\rm{amit}}}
\def\anit{{\rm{anit}}}
\def\arit{{\rm{arit}}}
\def\axit{{\rm{axit}}}
\def\sh{{\rm{sh}}}
\def\st{{\rm{st}}}
\def\swap{{\rm{swap}}}
\def\mt{{\frak{mt}}}
\def\ds{{\frak{ds}}}
\def\Q{{\Bbb Q}}
\def\u{{\bf u}}
\def\v{{\bf v}}
\def\w{{\bf w}}
\def\a{{\bf a}}
\def\b{{\bf b}}
\def\c{{\bf c}}

\centerline{{\bf Mould theory and the double shuffle Lie algebra structure}}
\bigskip
\centerline{Adriana Salerno\footnote{*}{The first author was partially supported by the NSF-AWM Mentoring Travel Grant.} and Leila Schneps}

\vskip 1cm
\noindent {\bf Abstract.} The real multiple zeta values 
$\zeta(k_1,\ldots,k_r)$ are known to form a
$\Q$-algebra; they satisfy a pair of well-known families of algebraic 
relations called the double shuffle relations.  In order to study the
algebraic properties of multiple zeta values, one can replace them by
formal symbols $Z(k_1,\ldots,k_r)$ subject only to the double shuffle
relations.  These form a graded Hopf algebra over $\Q$, and quotienting this 
algebra by products, one obtains a vector space.  A difficult theorem due to 
G. Racinet proves that this vector space carries the structure of a Lie 
coalgebra; in fact Racinet proved that the dual of this space is a Lie algebra, 
known as the double shuffle Lie algebra $\ds$.  

J. Ecalle developed a deep theory to explore combinatorial and algebraic
properties of the formal multiple zeta values. His theory is sketched out
in some publications (essentially [E1] and [E2]).  However, because of the depth 
and complexity of the theory, Ecalle did not include proofs of many of the most
important assertions, and indeed, even some interesting results are not
always stated explicitly.  The purpose of the present paper is to show how
Racinet's theorem follows in a simple and natural way from Ecalle's theory.  
This necessitates an introduction to the theory itself, which we have pared down 
to only the strictly necessary notions and results.
\vskip 1cm
\noindent {\bf \S 1. Introduction}
\vskip .3cm
In his doctoral thesis from 2000, Georges Racinet ([R1], see also [R2]) proved 
a remarkable theorem using astute combinatorial and algebraic 
reasoning. His proof was later somewhat simplified and streamlined by H. Furusho
([F]), but it remains really difficult to grasp the essential key that makes
it work.  The purpose of this article is to show how Ecalle's theory
of moulds yields a very different and deeply natural proof of the same result.
The only difficulty is to enter into the universe of moulds and learn its
language; the theory is equipped with a sort of standard all-purpose ``toolbox''
of objects and identities which, once acquired, serve to prove all kinds of
results, in particular the one we consider in this paper. Therefore, the goal of this article is not only to present the mould-theoretic proof of Racinet's theorem, but also to provide an initiation into mould theory in general.
Ecalle's seminal article on the subject is [E1], and a detailed introduction
with complete proofs can be found in [S]; the latter text will be referred
to here for some basic lemmas.
\vskip .2cm
We begin by recalling the definitions necessary to state Racinet's theorem.
\vskip .2cm
\noindent {\bf Definition.} Let $u,v$ be two monomials in $x$ and $y$.  Then
the commutative {\it shuffle product} $\sh(u,v)$ is defined recursively by 
$\sh(u,v)=u$ if $v=1$ and $v$ if $u=1$; otherwise, writing
$u=Xu'$ and $v=Yv'$ where $X,Y\in \{x,y\}$ represents the first letter
of the word, we have the recursive rule
$$\sh(Xu,Yv)=\{X\,\sh(u,Yv)\}\cup \{Y\,\sh(Xu,v)\}.\eqno(1.1)$$
If $u,v$ are two words ending in $y$, we can write them uniquely as words in the
letters $y_i=x^{i-1}y$.  The {\it stuffle product} of $u,v$ is defined by
$\st(u,v)=u$ if $v=1$ and $v$ if $u=1$, and
$$\st(y_iu,y_jv)=\{y_i\,\st(u,y_jv)\}\cup \{y_j\,\st(y_iu,v)\}\cup \{y_{i+j}\,\st(u,v)\},\eqno(1.2)$$
where $y_i$ and $y_j$ are respectively the first letters of the words $u$ and
$v$ written in the $y_j$.
\vskip .2cm
\noindent {\bf Definition.} The {\it double shuffle space} $\ds$ is the space
of polynomials $f\in {\Bbb Q}\langle x,y\rangle$, the polynomial ring on
two non-commutative variables $x$ and $y$, of degree $\ge 3$ and satisfying the
following two properties:
\vskip .2cm
\noindent (1) The coefficients of $f$ satisfy the {\it shuffle relations}
$$\sum_{w\in sh(u,v)} (f|w)=0,\eqno(1.3)$$
where $u,v$ are words in $x,y$ and $sh(u,v)$ is the set of words obtained
by shuffling them.  This condition is equivalent to the assertion that
$f$ lies in the free Lie algebra ${\rm Lie}[x,y]$, a fact that is easy
to see by using the characterization of Lie polynomials in the non-commutative
polynomial ring $\Q\langle x,y\rangle$ as those that are ``Lie-like'' under
the coproduct $\Delta$ defined by $\Delta(x)=x\otimes 1+1\otimes x$
and  $\Delta(y)=y\otimes 1+1\otimes y$, i.e. such that 
$\Delta(f)=f\otimes 1+1\otimes f$ (cf. [Se, Ch. III, Thm. 5.4]).  Indeed, 
when the property of being Lie-like
under $\Delta$ is expressed explicitly on the coefficients of $f$ it is
nothing other than the shuffle relations (1.3).
\vskip .2cm
\noindent (2) Let $f_*=\pi_y(f)+f_{\rm corr}$, where $\pi_y(f)$ is the
projection of $f$ onto just the words ending in $y$, and
$$f_{\rm corr}=\sum_{n\ge 1} {{(-1)^{n-1}}\over{n}}(f|x^{n-1}y)y^n.\eqno(1.4)$$
Considering $f_*$ as being rewritten in the variables $y_i=x^{i-1}y$, the coefficients of $f_*$ satisfy the {\it stuffle relations}:
$$\sum_{w\in st(u,v)} (f_*|w)=0,\eqno(1.5)$$
where $u$ and $v$ are words in the $y_i$.
\vskip .3cm
For every $f\in {\rm Lie}[x,y]$, define a derivation $D_f$ of
${\rm Lie}[x,y]$ by setting it to be
$$D_f(x)=0,\ \ D_f(y)=[y,f]$$
on the generators.  Define the {\it Poisson (or Ihara) bracket} on (the underlying
vector space of) ${\rm Lie}[x,y]$ by
$$\{f,g\}=[f,g]+D_f(g)-D_g(f).\eqno(1.6)$$
This definition corresponds naturally to the Lie bracket on the space of
derivations of ${\rm Lie}[x,y]$; indeed, it is easy to check that
$$[D_f,D_g]=D_f\circ D_g-D_g\circ D_f=D_{\{f,g\}}.\eqno(1.7)$$
\vskip .1cm
\noindent {\bf Racinet's Theorem.} {\it The double shuffle space $\ds$ is
a Lie algebra under the Poisson bracket.}
\vskip .3cm
The goal of this paper is to give the mould-theoretic proof of this result,
which first necessitates rephrasing the relevant definitions in terms
of moulds.  The paper is organized as follows. In \S 2, we give 
basic definitions from mould theory that will be used throughout the rest of 
the paper, and in \S 3 we define dimorphy and consider the main dimorphic 
subspaces related to double shuffle.  In \S 4 we give
the dictionary between mould theory and the double shuffle situation.
In \S 5 we give some of the definitions and basic results on the group aspect
of mould theory. In \S 6 we describe the special mould $pal$ that lies at 
the heart of much of mould theory, and introduce Ecalle's fundamental 
identity.  The final section \S 7 contains the simple and elegant proof of 
the mould version of Racinet's theorem.  Sections \S\S 2, 3, 5 and 6 can 
serve as a short introduction to the basics of mould theory; a much more
complete version with full proofs and details is given in [S], which
is cited for some results.  Every mould-theory definition in 
this paper is due to Ecalle, as are all of the statements, although some of
these are not made explicitly in his papers, but used as assumptions. 
Our contribution has been firstly to provide complete proofs of many 
statements which are either nowhere proved in his articles or proved by 
arguments that are difficult to understand (at least by us), secondly
to pick a path through the dense forest of his results that leads most
directly to the desired theorem, and thirdly, to give the dictionary that
identifies the final result with Racinet's theorem above.

In order to preserve the expository flow leading to the proof of the main 
theorem, we have chosen consign the longer and more technical proofs
to appendices or, for those that already appear in [S], to simply give the 
reference.
\vskip .5cm
\noindent {\bf \S 2. Definitions for mould theory}
\vskip .5cm
This section constitutes what could be called the ``first drawer'' of the mould
toolbox, with only the essential definitions of moulds, some operators on
moulds, and some mould symmetries.  We work over a base field $K$, and let 
$u_1,u_2,\ldots$ be a countable set of indeterminates, and $v_1,v_2,\ldots$ 
another. 
\vskip .3cm
\noindent {\it Moulds.} A {\it mould} in the variables $u_i$ is a family 
$A=(A_r)_{r\ge 0}$ of functions of the $u_i$, where each $A_r$ is a 
function of $u_1,\ldots,u_r$.  We call $A_r$ the {\it depth $r$} 
component of the mould. In this paper we let $K=\Q$, and in fact we 
consider only rational-function valued moulds, i.e. we have 
$A_r(u_1,\ldots,u_r)\in \Q(u_1,\ldots,u_r)$ for $r\ge 0$.
Note that $A_0(\emptyset)$ is a constant. We often drop the index
$r$ when the context is clear, and write $A(u_1,\ldots,u_r)$.
Moulds can be added and multiplied by scalars componentwise, so the set of
moulds forms a vector space. A mould in the $v_i$ is defined identically
for the variables $v_i$.
\vskip .3cm
Let $\ARI$ (resp. $\overline{\ARI}$) denote the space of moulds in the $u_i$ 
(resp. in the $v_i$) such that $A_0(\emptyset)=0$. These two vector spaces are obviously
isomorphic, but they will be equipped with very different Lie algebra
structures.  We use superscripts on $\ARI$ to denote the type of moulds
we are dealing with; in particular $\ARI^{pol}$ denotes the space of 
polynomial-valued moulds, and $\ARI^{rat}$ denotes the space of 
rational-function moulds.  
\vskip .5cm
\noindent {\it Operators on moulds.} We will use the following operators on
moulds in $\ARI$:
\vskip .2cm
$\neg(A)(u_1,\ldots,u_r)=A(-u_1,\ldots,-u_r)$\hfill (2.1)
\vskip .2cm
$\push(A)(u_1,\ldots,u_r)=A(-u_1-\cdots-u_r,u_1,\ldots,u_{r-1})$\hfill(2.2)
\vskip .2cm
$\mantar(A)(u_1,\ldots,u_r)=(-1)^{r-1}A(u_r,\ldots,u_1)$.\hfill(2.3)
\vskip .2cm
We also introduce the swap, which is a map from $ARI$ to $\overline{ARI}$
given by
\vskip .2cm
$\swap(A)(u_1,\ldots,u_r)=A(v_r,v_{r-1}-v_r,v_{r-2}-v_{r-1},\ldots,v_1-v_2)$,
\hfill(2.4)
\vskip .2cm
\noindent and its inverse, also called $\swap$, from $\overline{ARI}$ to
$ARI$:
\vskip .2cm
$\swap(A)(v_1,\ldots,v_r)=A(u_1+\cdots+u_r,u_1+\cdots+u_{r-1},\ldots,u_1+u_2,u_1).
$\hfill(2.5)
\vskip .3cm
Thanks to this formulation, which is not ambiguous since to know which swap
is being used it suffices to check whether swap is being applied to a mould
in $ARI$ or one in $\overline{ARI}$, we can treat swap like an involution:
$\swap\circ \swap={\rm id}$.
\vskip .3cm
Let us now introduce some notation necessary for the Lie algebra structures on $\ARI$ and
$\overline{\ARI}$. 
\vskip .5cm
\noindent {\it Flexions.}  Let $\w=(u_1,\cdots,u_r)$.
For every possible way of cutting the word $\w$ into three (possibly empty) subwords
$\w=\a\b\c$ with
$$\a=(u_1,\ldots,u_k),\ \ 
\b=(u_{k+1},\ldots,u_{k+l}),\ \
\c=(u_{k+l+1},\ldots,u_r),$$
set
$$\cases{\a\rceil=(u_1,u_2,\cdots,u_k+u_{k+1}+\cdots+u_{k+l})
&if $\b\ne\emptyset$, otherwise $\a\rceil=\a$\cr
\lceil\c=(u_{k+1}+\cdots+u_{k+l+1},u_{k+l+2},\cdots,u_r)
&if $\b\ne \emptyset$, otherwise $\lceil\c=\c$.}$$
If now $\w=(v_1,\ldots,v_r)$ is a word in the $v_i$, then for every 
decomposition $\w=\a\b\c$ with
$$\a=(v_1,\ldots,v_k),\ \ 
\b=(v_{k+1},\ldots,v_{k+l}),\ \
\c=(v_{k+l+1},\ldots,v_r),$$
we set
$$\cases{
\b\rfloor=(v_{k+1}-v_{k+l+1},v_{k+2}-v_{k+l+1},\ldots,v_{k+l}-v_{k+l+1})&if $\c\ne\emptyset$, otherwise $\b\rfloor=\b$\cr
\lfloor\b=(v_{k+1}-v_k,v_{k+2}-v_k,\ldots,v_{k+l}-v_k)&if $\a\ne\emptyset$,
otherwise $\lfloor\b=\b$.}$$
\vfill\eject
\noindent {\it Operators on pairs of moulds.} 
For $A,B\in \ARI$ or $A,B\in\overline{\ARI}$, we set
\vskip .2cm
$\mu(A,B)(\w)=\displaystyle\sum_{\w=\a\b} A(\a)B(\b)-B(\a)A(\b)$\hfill(2.6)
\vskip .2cm
$\lu(A,B)=\mu(A,B)-\mu(B,A)$.\hfill(2.7)
\vskip .3cm
For any mould $B\in\ARI$, we
define two operators on $\ARI$, $\amit(B)$ and $\anit(B)$, defined by
\vskip .3cm
$\bigl({\amit}(B)\cdot A\bigr)(\w)=
\displaystyle{\sum_{{{\w=\a\b\c}\atop{\b,\c\ne \emptyset}}}}
A(\a\lceil\c)B(\b)$

$\bigl({\anit}(B)\cdot A\bigr)(\w)=
\displaystyle{\sum_{{{\w=\a\b\c}\atop{\a,\b\ne \emptyset}}}}
A(\a\rceil\c)B(\b)$.\hfill(2.8)

The operators $\amit$ and $\anit$ are derivations of $\ARI$ for the
$\lu$-bracket (see [S], Prop.~2.2.1).  We define a third derivation,
$\arit(B)$, by
$$\bigl({\arit}(B)\cdot A\bigr)(\w)=\amit(B)\cdot A-\anit(B)\cdot A.\eqno(2.9)$$
\vskip .2cm
\noindent If $B\in\overline{\ARI}$ we have derivations of $\overline{\ARI}$
given by
\vskip .3cm
$\bigl({\amit}(B)\cdot A\bigr)(\w)=
\displaystyle{\sum_{{{\w=\a\b\c}\atop{\b,\c\ne \emptyset}}}}
A(\a\c)B(\b\rfloor)$

$\bigl({\anit}(B)\cdot A\bigr)(\w)=
\displaystyle{\sum_{{{\w=\a\b\c}\atop{\a,\b\ne \emptyset}}}}
A(\a\c)B(\lfloor\b)$,

\noindent and again we define the derivation $\arit(B)$ as in (2.9).
\vskip .3cm
Finally, we set 
\vskip .2cm
$\ari(A,B)=\arit(B)\cdot A+\lu(A,B)-\arit(A)\cdot B$.\hfill(2.10)
\vskip .3cm
Since $\arit$ is a derivation for $\lu$, the $\ari$-operator is easily 
shown to be a Lie bracket.  Note that although we use the same notation 
$\ari$ for the Lie brackets on both $\ARI$ and $\overline{\ARI}$, they 
are two different Lie brackets on two different spaces. Indeed, while 
some formulas and properties (such as $\mu$, or alternality, see (2.11) 
below) are written identically for $\ARI$ and $\overline{\ARI}$, others, in 
particular all those that use flexions, are very different. 
\vskip .5cm
\noindent {\it Symmetries.} A mould in $\ARI$ (resp. $\overline{\ARI}$) is said to be
{\it alternal} if for all words $\u,\v$ in the $u_i$ (resp. $v_i$),
$$\sum_{\w\in \sh(\u,\v)} A(\w)=0.\eqno(2.11)$$
The relations in (2.11) are known as the {\it alternality} relations, and
they are identical for moulds in $\ARI$ and $\overline{\ARI}$.  Let us now 
define the {\it alternility relations}, which are only applicable to moulds 
in $\overline{\ARI}$.  Just as the alternality conditions are the mould 
equivalent of the shuffle relations, the alternility conditions are the 
mould equivalent of the stuffle relations, translated in 
terms of the alphabet $\{v_1,v_2,\ldots\}$ as follows. Let
$Y_1=(y_{i_1},\ldots,y_{i_r})$ and $Y_2=(y_{j_1},\ldots,y_{j_s})$ be two 
sequences; for example, we consider $Y_1=(y_i,y_j)$ and $Y_2=(y_k,y_l)$.
Let $w$ be a word in the stuffle product $st\bigl(Y_1,Y_2\bigr)$, which
in our example is the 13-element set
$$\bigl\{(y_i,y_j,y_k,y_l),(y_i,y_k,y_j,y_l),(y_i,y_k,y_l,y_j),(y_k,y_i,y_j,y_l),
(y_k,y_i,y_l,y_j),(y_k,y_l,y_i,y_j),$$
$$(y_i,y_{j+k},y_l), (y_{i+k},y_j,y_l), (y_i,y_k,y_{j+l}),
(y_{i+k},y_l,y_j), (y_k,y_i,y_{j+l}), (y_k,y_{i+l},y_j),$$
$$(y_{i+k},y_{j+l})\bigr\}.\eqno(2.12)$$
To each such word we associate an alternility term for the mould $A$,
given by associating the tuple $(v_1,v_2,v_3,v_4)$ to
the ordered tuple $(y_i,y_j,y_k,y_l)$ and taking
$${{1}\over{(v_i-v_j)}}\bigl(A(\ldots,v_i,\ldots)-A(\ldots,v_j,\ldots)\bigr)
\eqno(2.13)$$
for each contraction occurring in the word $w$. For instance in our example
we have the six alternility terms
$$A(v_1,v_2,v_3,v_4),A(v_1,v_3,v_2,v_4),A(v_1,v_3,v_4,v_2),
A(v_3,v_1,v_2,v_4),$$ $$A(v_3,v_1,v_4,v_2),A(v_3,v_4,v_1,v_2)\eqno(2.14)$$
corresponding to the first six words in (2.12), the six terms
$${{1}\over{(v_2-v_3)}}\bigl(A(v_1,v_2,v_4)-A(v_1,v_3,v_4)\bigr),\ \ \ 
{{1}\over{(v_1-v_3)}}\bigl(A(v_1,v_2,v_4)-A(v_3,v_2,v_4)\bigr),$$
$${{1}\over{(v_2-v_4)}}\bigl(A(v_1,v_3,v_2)-A(v_1,v_3,v_4)\bigr),\ \ \ 
{{1}\over{(v_1-v_3)}}\bigl(A(v_1,v_4,v_2)-A(v_3,v_4,v_2)\bigr),$$
$${{1}\over{(v_2-v_4)}}\bigl(A(v_3,v_1,v_2)-A(v_3,v_1,v_4)\bigr),\ \ \
{{1}\over{(v_1-v_4)}}\bigl(A(v_3,v_1,v_2)-A(v_3,v_4,v_2)\bigr)\eqno(2.15)$$
corresponding to the next six words, and the final term
$${{1}\over{(v_1-v_3)(v_2-v_4)}}\bigl(A(v_1,v_2)-A(v_3,v_2)-A(v_1,v_4)
+A(v_3,v_4)\bigr)\eqno(2.16)$$
corresponding to the final word with the double contraction.
Writing $A_{\w}$ for the alternility term of $A$ associated to a word $\w$
in the stuffle product $\st(Y_1,Y_2)$, the {\it alternility relation} 
associated to the pair $(Y_1,Y_2)$ on $A$ is given by 
$$\sum_{\w\in \st(Y_1,Y_2)} A_{\w}=0.\eqno(2.17)$$  
Let $A_{r,s}$ denote the left-hand side of (2.17).  Note that indeed,
$A_{r,s}$ does not depend on the actual sequences $Y_1$ and $Y_2$, but 
merely on the number of letters in $Y_1$ and in $Y_2$. For example
when $r=s=2$, the alternility sum $A_{2,2}$ is given by the sum of
the terms (2.14)-(2.16) above.  Furthermore, like 
for the shuffle, we may assume that $r\le s$ by symmetry.  Thus we have the 
following definition: a mould in $\overline{\ARI}$ is said to be 
{\it alternil} if it satisfies the alternility relation $A_{r,s}=0$
for all pairs of integers $1\le r\le s$.
\vskip .3cm
\noindent {\bf Remark.} If $A$ is a polynomial-valued mould, then the
alternility terms are all polynomial.
To see this, it suffices to note that setting $v_i=v_j$ in (2.13) gives a
zero in the numerator that cancels out the pole in the denominator.
\vfill\eject
\noindent {\bf \S 3. Lie subalgebras of \ARI}
\vskip .3cm
In this section, we show that the spaces of moulds satisfying certain important 
symmetry properties are closed under the $\ari$-bracket. 
\vskip .3cm
\noindent {\bf Definitions.} Let $\ARI_{al}$ denote the set of alternal moulds.
Let $\ARI_{al/al}$ (resp. $\ARI_{al/il}$) denote the set of alternal moulds with 
alternal (resp. alternil) swap.  Let
$\ARI_{al*al}$ (resp. $\ARI_{al*il}$) denote the set of alternal moulds
whose swap is alternal (resp. alternil) up to addition of a constant-valued
 mould.
Finally, let $\ARI_{\underline{al}/\underline{al}}$ 
(resp. $\ARI_{\underline{al}*\underline{al}}$, $\ARI_{\underline{al}/\underline{il}}$, $\ARI_{\underline{al}*\underline{il}}$) denote the subspace of $\ARI_{al/al}$
(resp. $\ARI_{al*al}$, $\ARI_{al/il}$, $\ARI_{al*il}$) consisting of moulds $A$ such 
that $A_1$ is an even function, i.e.  $A(-u_1)=A(u_1)$.  
\vskip .3cm
The first main theorem of this paper is the following result, which is
used constantly in Ecalle's work although no explicit proof appears to have
been written down, and the proof is by no means as easy as one might imagine.
\vskip .3cm
\noindent {\bf Theorem 3.1.} {\it The subspace $\ARI_{al}\subset \ARI$ of alternal moulds forms a Lie algebra under the $\ari$-bracket, as does the subspace
$\overline{\ARI}_{al}$ of $\overline{\ARI}$.}
\vskip .3cm
\noindent The full proof is given in Appendix A. The idea is as follows:
if $C=\ari(A,B)$, then by (2.10) it is enough to show separately that if 
$A$ and $B$ are alternal then $\lu(A,B)$ is alternal and $\arit(B)\cdot A$ 
is alternal. This is done via a combinatorial manipulation that is 
fairly straightforward for $\lu$ but actually quite complicated for $\arit$.
\vskip .3cm
We next have a simple but important result on polynomial-valued moulds.
\vskip .3cm
\noindent {\bf Proposition 3.2.} {\it The subspace $ \ARI^{pol}$ of polynomial-valued moulds in $\ARI$ forms a Lie algebra under the $\ari$-bracket.}
\vskip .2cm
\noindent Proof. This follows immediately from the definitions of $\mu$, 
$\arit$ and $\ari$ in (2.6)-(2.9), as all the operations and flexions there
are polynomial.\hfill{$\diamondsuit$}
\vskip .3cm
Now we give another key theorem, the first main result concerning 
{\it dimorphy}, i.e. the double description of a mould by a symmetry property 
on it and another one on its swap. This result, again, is used repeatedly
by Ecalle but we were not able to find a complete proof in his papers, so we 
have reconstructed one here.
\vskip .3cm
\noindent {\bf Theorem 3.3.} {\it The subspaces $\ARI_{\underline{al}/\underline{al}}$ and $\ARI_{\underline{al}*\underline{al}}$ form Lie algebras under the $\ari$-bracket.}
\vskip .3cm
The proof is based on the following two propositions. 
\vskip .3cm
\noindent {\bf Proposition 3.4.} {\it If $A\in 
\ARI_{\underline{al}*\underline{al}}$,
then $A$ is $\neg$-invariant and $\push$-invariant.}
\vskip .3cm
The proof of this proposition is deferred to Appendix B.
\vskip .3cm
\noindent {\bf Proposition 3.5.} {\it If $A$ and $B$ are both $\push$-invariant moulds, then 
$$\swap\Bigl(\ari\bigl(\swap(A),\swap(B)\bigr)\Bigr)=\ari(A,B),\eqno(3.1)$$}
\par
\noindent Proof. Explicit computation using the flexions shows that
for all moulds $A,B\in \ARI$ we have the general formula: 
$$\swap\bigl(\ari(\swap(A),\swap(B))\bigr)=
\axit\bigl(B,-\push(B)\bigr)\cdot A-\axit\bigl(A,-\push(A)\bigr)\cdot B+\lu(A,B),\eqno(3.2)$$
where here $\ari$ is the Lie bracket on $\overline{\ARI}$, and 
$\axit$ is the operator on $\ARI$ defined for a general pair of moulds 
$B,C\in \ARI$ by the formula
$$\axit(B,C)\cdot A=\amit(B)\cdot A+\anit(C)\cdot A.$$
(See [S, \S 4.1] for complete details of this flexion computation.)
Comparing with (2.9) shows that $\arit(B)=\axit(B,-B)$. Thus if $A$ and 
$B$ are $\push$-invariant, (3.2) reduces to
$$\swap\Bigl(\ari\bigl(\swap(A),\swap(B)\bigr)\Bigr)=\arit(B)\cdot A-\arit(A)\cdot B+\lu(A,B),$$
which is exactly $\ari(A,B)$ by (2.10). 
\hfill{$\diamondsuit$}
\vskip .5cm
\noindent Proof of Theorem 3.3.
Using the propositions, the proof becomes reasonably easy.
We first consider the case where $A$, 
$B\in \ARI_{\underline{al}/\underline{al}}$. 
In particular $A$ and $B$ are alternal. Set $C=\ari(A,B)$. 
The mould $C$ is alternal by Theorem 3.1.  By Proposition 3.4, we know that
$A$ and $B$ are $\push$-invariant, so by Proposition 3.5 we have
$\swap(C)=\swap\bigl(\ari(A,B)\bigr)=\ari\bigl(\swap(A),\swap(B)\bigr)$. 
But this is also alternal by Theorem 3.1, so $C\in \ARI_{al/al}$. 
Furthermore, it follows directly from the defining formula for the 
$\ari$-bracket, which is additive in the mould depths, that if $C$ is an 
$\ari$-bracket of two moulds in $\ARI$, i.e. with constant term equal to 0, 
we must have $C(u_1)=0$, so $C\in \ARI_{\underline{al}/\underline{al}}$.

Now we consider the more general situation where $A,B\in 
\ARI_{\underline{al}*\underline{al}}$.  Let $A_0, B_0$ be the 
constant-valued moulds such that $\swap(A)+A_0$ and $\swap(B)+B_0$
are alternal.  From the definitions (2.6)-(2.9), we see that for any 
constant-valued mould $M_0$, we have $\arit(M_0)\cdot M=0$ and 
$\arit(M)\cdot M_0=\lu(M,M_0)$; thus by (2.10), we have $\ari(M,M_0)=0$.
Thus we find that 
$$\ari(A+A_0,B+B_0)=\ari(A,B)+\ari(A,B_0)+\ari(A_0,B)+\ari(A_0,B_0)=
\ari(A,B).\eqno(3.3)$$
Now, $A$ and $B$ are push-invariant by Proposition 3.4, and constant-valued
moulds are always push-invariant, so $A+A_0$ and $B+B_0$ are also
push-invariant; thus we have
$$\eqalign{\swap(C)&=\swap\bigl(\ari(A,B)\bigr)\cr
&=\swap\bigl(\ari(A+A_0,B+B_0)\bigr)\ \ \ {\rm by\ (3.3)}\cr
&=\ari\bigl(\swap(A+A_0),\swap(B+B_0)\bigr)\ \ \ {\rm by\ (3.1)}.}$$
But since $\swap$ preserves constant-valued moulds, we have
$\swap(A+A_0)=\swap(A)+A_0$ and $\swap(B+B_0)=\swap(B)+B_0$. These two
moulds are alternal by hypothesis, so by Theorem 3.1, their
$\ari$-bracket is alternal, i.e. $\swap(C)$ is alternal. Since as above 
we have $C(u_1)=0$, we find that in fact
$C$ is not just in $\ARI_{\underline{al}*\underline{al}}$
but in $\ARI_{\underline{al}/\underline{al}}$.  This completes the proof
of Theorem 3.3.
\hfill{$\diamondsuit$}
\vskip .3cm
We will see in the next section that the double shuffle space $\ds$ defined 
in \S 1 is isomorphic to the space of polynomial-valued moulds 
$\ARI^{pol}_{\underline{al}*\underline{il}}$, with the alternality
property translating shuffle and the alternility property translating
stuffle. Thus dimorphy is closely connected to double shuffle, but much 
more general, since the symmetry properties of alternality or alternility 
on itself or its swap can hold for any mould, not just polynomial ones.
\vskip .5cm
\noindent {\bf \S 4. Dictionary with the Lie algebra and
double shuffle framework}
\vskip .3cm
Let $C_i=ad(x)^{i-1}y\in{\Bbb Q}\langle x,y\rangle$.
By Lazard elimination, the subring ${\Bbb Q}\langle C_1,C_2,\ldots\rangle$,
which we denote simply by ${\Bbb Q}\langle C\rangle$,
is free on the $C_i$.  For $a\in {\Bbb Q}$, Let ${\Bbb Q}_a\langle C\rangle$ 
denote the subspace of polynomials in the $C_i$ with constant term equal
to $a$.  Define a linear map
$$\eqalign{ma:{\Bbb Q}_{0}\langle C\rangle&\buildrel\sim\over
\rightarrow \ARI^{pol}\cr
C_{a_1}\cdots C_{a_r}&\mapsto A_{a_1,\ldots,a_r}}\eqno(4.1)$$
where $A_{a_1,\ldots,a_r}$ is the polynomial mould concentrated in depth $r$
defined by
$$A_{a_1,\ldots,a_r}(u_1,\ldots,u_r)=(-1)^{a_1+\cdots+a_r-r}u_1^{a_1-1}\cdots u_r^{a_r-1}.\eqno(4.2)$$
This map $ma$ is trivially invertible and thus an isomorphism. Let
${\rm Lie}[C]$ denote the free Lie algebra ${\rm Lie}[C_1,C_2,\ldots]$ on
the $C_i$.  Note that, again by Lazard elimination, we can write 
${\rm Lie}[x,y]={\Bbb Q}x\oplus {\rm Lie}[C]$, and that
since the double shuffle space $\ds\subset {\rm Lie}[x,y]$ contains only
polynomials of degree $\ge 3$, we have 
$$\ds\subset {\rm Lie}[C]\subset {\Bbb Q}_0 \langle C\rangle.$$
\noindent {\bf Definition.} Let ${\cal{MT}}_0$ denote the Lie algebra 
whose underlying space is the space of polynomials
${\Bbb Q}_0\langle C\rangle$,
equipped with the Poisson bracket (1.6), and let $\mt$ denote the subspace
of Lie polynomials in the $C_i$, i.e. the vector space ${\rm Lie}[C]$
equipped with the Poisson bracket. Observe that $\mt$ is closed under
the Poisson bracket since if $f,g$ are Lie then so are $D_f(g)$, 
$D_g(f)$ and $[f,g]$, so $\mt$ is a Lie algebra.  The letters ``M-T'' stand for {\it twisted
Magnus} (cf. [R1]). Let ${\cal{MT}}$ denote the universal enveloping 
algebra of $\mt$, isomorphic as a vector space to 
${\Bbb Q}\langle C\rangle$,
but equipped with a multiplication coming from the pre-Lie law 
$f\odot g=fg-D_g(f)$ on $\mt$. We also define the twisted Magnus group as
the exponential $MT=\exp^\odot(\mt)$, where
$exp^\odot(f)=\sum_{n\ge 0} {{1}\over{n!}} f^{\odot n}$. 
\vskip .3cm
\noindent {\bf Theorem 4.1.} (Racinet) {\it The linear isomorphism 
(4.1) is a Lie algebra isomorphism 
$$ma:{\cal{MT}}_0\buildrel\sim\over\rightarrow\ARI^{pol},\eqno(4.3)$$
and it restricts to a Lie algebra isomorphism of the Lie subalgebras
$$ma:{\mt}\buildrel\sim\over\rightarrow\ARI^{pol}_{al}.\eqno(4.4)$$}\par
\noindent Proof.  In view of the fact that
$ma$ is invertible as a linear map, the isomorphism (4.3) follows 
from the following identity relating the Poisson bracket and the 
$\ari$-bracket on polynomial-valued moulds, which was proven by Racinet 
in his thesis ([R1, Appendix A], see also [S, Corollary 3.3.4]):
$$ma\bigl(\{f,g\}\bigr)=\ari\bigl(ma(f),ma(g)\bigr),\eqno(4.5)$$
The isomorphism (4.4), identifying Lie polynomials with
alternal polynomial moulds, follows from a standard argument that 
we indicate briefly, as it is merely an adaptation to ${\rm Lie}[C]$
of the similar argument following the definition of the shuffle relations 
in (1.3). Let $\Delta$ denote 
the standard cobracket on ${\Bbb Q}\langle C\rangle$ defined by $\Delta(C_i)=C_i\otimes 1+1\otimes C_i$.  Then 
the Lie subspace ${\rm Lie}[C]$ of the polynomial algebra
${\Bbb Q}\langle C\rangle$ is the space of primitive elements for 
$\Delta$, i.e. elements $f\in {\rm Lie}[C]$ satisfying 
$\Delta(f)=f\otimes 1+1\otimes f$.  This condition on $f$ is given explicitly
on the coefficients of $f$ by the family of shuffle relations 
$$\sum_{D\in sh(C_{a_1}\cdots C_{a_r},C_{b_1}\cdots C_{b_s})} (f|D)=0,$$
where $(f|D)$ denotes the coefficient in the polynomial $f$ of the
monomial $D$ in the $C_i$.  But these conditions are exactly equivalent to 
the alternality relations 
$$\sum_{D\in sh((a_1,\ldots,a_r),(b_1,\ldots,b_s))} ma(f)(D)=0,$$
proving (4.4).\hfill{$\diamondsuit$}
\vskip .3cm
\noindent {\bf Theorem 4.3.} {\it The linear isomorphism (4.4) restricts
to a linear isomorphism of the subspaces
$$ma:\ds\buildrel\sim\over\rightarrow \ARI^{pol}_{\underline{al}*\underline{il}}.\eqno(4.6)$$}
\noindent Proof.  By (4.4), since $\ds\subset \mt$, we have $ma:\ds
\hookrightarrow \ARI^{pol}_{al}$.  If an element 
$f\in\ds$ has a depth 1 component, i.e. if the coefficient of $x^{n-1}y$ 
in $f$ is non-zero, then $n$ is odd; this is a simple consequence of 
solving the depth 2 stuffle relations (see [C, Theorem 2.30 (i)] for 
details). Thus, if the mould $ma(f)$ has a depth 1 component,
it will be an even function, since by the definition of $ma$ the degree of 
$ma(f)(u_1)$ is equal to the degree of $f$ minus 1.  This shows that
$ma$ maps $\ds$ to moulds that are even in depth 1, i.e.
$$ma:\ds\hookrightarrow\ARI^{pol}_{\underline{al}}.$$
It remains only to show that if $f\in \ds$ then $\swap\bigl(ma(f)\bigr)$
is alternil up to addition of a constant mould, i.e. that the stuffle
conditions (1.5) imply the alternility of 
$\swap\bigl(ma(f)\bigr)$.  

By additivity, we may assume that $f$ is of homogeneous degree $n$. 
Let $C$ be the constant mould concentrated in depth $n$ given by 
$C(u_1,\ldots,u_n)={{(-1)^{n-1}}\over{n}}(f|x^{n-1}y)$, and let
$A=mi(f)+C$.  Then $A$ is concentrated in depths $\le n$, and $A$ is
obtained directly from the polynomial $f_*$ by replacing
monomials in the $y_i$ by corresponding monomials in the $v_i$: explicitly,
if the depth $r$ part of $f_*$ is given by
$$(f_*)^r=\sum_{{\bf a}=(a_1,\ldots,a_r)} c_{{\bf a}}\,y_{a_1}\cdots y_{a_r},
\eqno(4.7)$$
then the depth $r$ part of $A$ is given by
$$A(v_1,\ldots,v_r)=\sum_{{\bf a}=(a_1,\ldots,a_r)} c_{{\bf a}}\,v_1^{a_1-1}\cdots v_r^{a_r-1}\eqno(4.8)$$
Tis follows from the definition of $mi(f)$ as $swap\bigl(ma(f)\bigr)$
(see [R1, Appendix A] or [S, (3.2.6)] for details). 
So we need to show that the stuffle relations (1.5) on $f_*$ are equivalent
to the alternility of $A$.

For any pair of integers $1\le r\le s$, let $A_{r,s}$ denote the 
alternility sum associated to the mould $A$ as in (2.17). By definition,
$A$ is alternil if and only if $A_{r,s}=0$ for all 
pairs $1\le r\le s$. Recall from \S 2 that the alternility sum
$A_{r,s}$ is a polynomial in $v_1,\ldots,v_{r+s}$ obtained by summing up 
polynomial terms in one-to-one correspondence with the terms of the
stuffle of two sequences of lengths $r$ and $s$.
By construction, the coefficient of a monomial $w=v_1^{b_1-1}
\cdots v_{r+s}^{b_{r+s}-1}$ in the alternility term corresponding to 
to a given stuffle term is equal to the coefficient in $f_*$ of the
stuffle term itself.  This follows from a direct calculation obtained by
expanding the alternility terms; for example, the alternility term
corresponding to the stuffle term $(y_i,y_{j+k},y_l)$ in (2.12)
is given by
$${{1}\over{v_2-v_2}}\bigl(A(v_1,v_2,v_4)-A(v_1,v_3,v_4)\bigr)$$
(see (2.14)), whose polynomial expansion is given by
$$\sum_{{\bf a}=(a_1,a_2,a_3)} c_{{\bf a}}v_1^{a_1-1}\bigl(\sum_{m=0}^{a_2-2}
v_2^mv_3^{a_2-2-m}\bigr)v_4^{a_3-1},$$
and the coefficient of the monomial 
$v_1^{i-1}v_2^{j-1}v_3^{k-1}v_4^{l-1}$ in this alternility term corresponds
to $a_1-1=i-1$, $m=j-1$, $a_2-2-m=k-1$ and $a_3-1=l-1$, i.e.
$a_1=i$, $a_2=j+k$, $a_3=l$, so it is equal to $c_{i,j+k,l}$ which is
exactly the coefficient $(f_*|y_iy_{j+k}y_l)$ in (4.7).
The alternility sum is equal to zero if and only the coefficient of each
monomial in $v_1,\ldots,v_{r+s}$ is equal to zero, which is thus equivalent
to the full set of stuffle relations on $f_*$.
\hfill{$\diamondsuit$} 
\vskip .3cm
In view of (4.5) and (4.6), a mould-theoretic proof of Racinet's theorem 
consists in proving that $\ARI^{pol}_{\underline{al}*\underline{il}}$ is a 
Lie algebra under the $\ari$-bracket.  To prove this mould-theoretic
version, we need to make use of the Lie group $\GARI$ associated to
$\ARI$, defined in the next section.  In \S 6 we give the necessary results 
from Ecalle's theory, and the theorem is proved in \S 7.
\vfill\eject
\noindent {\bf \S 5. The group {\it \GARI}}
\vskip .3cm
In this section we introduce several notions on the group $\GARI$ of moulds 
with constant term 1, which are group analogs of the Lie notions introduced
in \S 2.  To move from the Lie algebra ARI to the associated
group GARI, Ecalle introduces a {\it pre-Lie} law on ARI,
defined as follows:
$$\preari(A,B) = \arit(B)\cdot A+\mu(A,B),\eqno(5.1)$$
where $\arit$ and $\mu$ are as defined in (2.9) and (2.6). Using these,
he defines an exponential map on ARI in the standard way:
$$\expari(A) = \sum_{n\ge 0} {{1}\over{n!}}\,\preari(\underbrace{A,\ldots,A}_n),\eqno(5.2)$$
where 
$$\preari(\underbrace{A,\ldots,A}_n)=\preari(\preari(\underbrace{A,\ldots,A}_{n-1}),A).$$
This map is the exponential isomorphism $\expari:\ARI\rightarrow \GARI$, 
where $\GARI$ is nothing other than the group of all moulds with 
constant term equal to $1$, equipped with the multiplication law, denoted
gari, that comes as always from the Campbell-Hausdorff law ${\rm ch}(\cdot ,\cdot)$
on $\ARI$:
$$\gari\bigl(\expari(A),\expari(B)\bigr)=\expari({\rm ch}(A,B)).\eqno(5.3)$$

The gari-inverse of a mould $B\in \GARI$ is denoted $\invgari(B)$.  The inverse isomorphism of $\expari$ is denoted by $\logari$.  
\vskip .2cm
Like all Lie algebras, $\ARI$ is equipped with an action of the associated
group $\GARI$, namely the standard adjoint action, denoted
$\adari$ (Ecalle denotes it simply adari, but we have modified it to
stress the fact that it is represents the adjoint action of the group
$GARI$ on $ARI$):
$$\eqalign{\adari(A)\cdot B&=\gari\bigl(\preari(A,B),\invgari(A)\bigr)\cr
&= {{d}\over{dt}}|_{t=0}\ \gari\bigl(A,\expari(tB),\invgari(A)\bigr)\cr
&=B+\ari\bigl(\logari(A),B\bigr)+{{1}\over{2}}\ari\bigl(\logari(A),\ari\bigl(\logari(A),B\bigr)+\cdots}\eqno(5.4)$$

Finally, to any mould $A\in \GARI$ (i.e.~any mould in the $u_i$ with 
constant term 1), Ecalle associates an automorphism $\ganit(A)$ of the ring 
of all moulds in the $u_i$ under the $\mu$-multiplication which is
just the exponential of the derivation $\anit\bigl(\logari(A)\bigr)$.

%$$\bigl(\ganit(A)\cdot B\bigr)({\bf w})=\sum B({\bf b}^1\rceil\cdots
%{\bf b}^s\rceil)A({\bf c}^1)\cdots A({\bf c}^s),\eqno(5.5)$$
%where the sum runs over the decompositions of the word 
%${\bf w}=(u_1,\ldots,u_r)$ ($r\ge 1$) as 
%$${\bf w}={\bf b}^1{\bf c}^1\cdots {\bf b}^s{\bf c}^s,\ \ (s\ge 1)$$
%where all ${\bf b}^i$ and ${\bf c}^i$ are non-empty words except
%possibly for ${\bf c}^s$.
\vskip .2cm
The analogous objects exist for moulds in the $v_i$. If $\preari$ denotes 
the pre-Lie law on $\overline{\ARI}$ given by (5.1) (but for the
derivation $\arit$ of $\ARI$), then the formula (5.2) defines an analogous 
exponential isomorphism $\overline{\ARI}\rightarrow \GVARI$, where $\GVARI$ 
consists of all moulds in the variables $v_i$ with constant term 1 and 
multiplication determined by (5.3) (note that this definition depends on 
that of $\arit$, so just as the Lie bracket $\ari$ is different for
$\ARI$ and $\overline{\ARI}$, the multiplication is different for 
$\GARI$ and $\GVARI$).  As above, we let the automorphism $\ganit(A)$ 
of $\GVARI$ associated to each $A\in \GVARI$ be defined as the exponential 
of the derivation $\anit\bigl(\logari(A)\bigr)$ of $\overline{\ARI}$.
%defining formula for $\ganit(A)$ in the $v_i$ is
%$$\bigl(\ganit(A)\cdot B\bigr)({\bf w})=\sum B({\bf b}^1\cdots
%{\bf b}^s)A(\lfloor{\bf c}^1)\cdots A(\lfloor{\bf c}^s),\eqno(5.6)$$
\vskip .2cm
\noindent {\bf Definition.} 
A mould $A\in\GARI$ (resp. $\GVARI$) is {\it symmetral} if for all words $\u,\v$ 
in the $u_i$ (resp. in the $v_i$), we have
$$\sum_{\w\in\sh(\u,\v)} A(\w)=A(\u)A(\v).\eqno(5.7)$$
Following Ecalle, we write $\GARI_{as}$ (resp. $\GVARI_{as}$) for the
set of symmetral moulds in $\GARI$ (resp. $\GVARI$).  The property of 
{\it symmetrality} is the group equivalent of alternality; in particular, 
$$A\in \ARI_{al}\ {\rm (resp.}\ \overline{\ARI}_{al}{\rm )}\ \Leftrightarrow
\ \expari(A)\in \GARI_{as}\  {\rm (resp.}\ \GVARI_{as}{\rm )}.
\eqno(5.8)$$
\vskip .5cm
\noindent {\bf Remark.} Let $MT$ denote the {\it twisted Magnus group} of
power series in ${\Bbb Q}\langle\langle C_1,C_2,\ldots\rangle\rangle$ with
constant term 1, identified with the exponential of the twisted Magnus
Lie algebra $\mt$ defined by
$$\exp^\odot(f)=\sum_{n\ge 0} {{1}\over{n!}}\,f^{\odot n}$$
for $f\in\mt$, where $\odot$ is the pre-Lie law 
$$f\odot g=fg+D_f(g)\eqno(5.9)$$
defined for $f,g\in \mt$ (see \S 4).  
The group $MT$ is equipped with the twisted Magnus multiplication
$$\bigl(f\odot g\bigr)(x,y)=f(x,gyg^{-1})g(x,y).\eqno(5.10)$$
Notice that it makes sense to use the same symbol $\odot$ for (5.9) and
(5.10), because in fact $\odot$ is the multiplication on the completion
of the universal enveloping algebra of $\mt$, and (5.9) and (5.10) merely
represent the particular expressions that it takes on two elements of
$\mt$ resp. two elements of $MT$.

The multiplication (5.10) corresponds to the $\gari$-multiplication in the sense
that the map $ma$ defined in (4.1) yields a group isomorphism 
$MT\buildrel\sim\over\rightarrow\GARI^{pol}$.  If $g\in MT$, then the
automorphism $\ganit\bigl(ma(g)\bigr)$ is the $\GARI$-version of the automorphism of
$MT$ given by mapping $x\mapsto x$ and $y\mapsto yg$.  
\vskip .2cm
The fact of having
non-polynomial moulds in $\GARI$ gives enormously useful possibilities 
of expanding the familiar symmetries and operations (derivations, shuffle
and stuffle relations etc.) to a broader situation. 
In particular, the next section contains some of Ecalle's most important 
results in mould multizeta theory, which make use of moulds with 
denominators and have no analog within the usual polynomial framework.
\vfill\eject
\noindent {\bf \S 6.  The mould pair $pal/pil$ and Ecalle's fundamental 
identity} 
\vskip .3cm
In this section we enter into the ``second drawer'' of Ecalle's powerful
toolbox, with the mould pair $pal/pil$ and Ecalle's fundamental identity. 
\vskip .3cm
\noindent {\bf Definition.}  Let $dupal$ be the mould defined explicitly
by the following formulas: $dupal(\emptyset)=0$ and for $r\ge 1$,
$$dupal(u_1,\ldots,u_r)={{B_r}\over{r!}}
{{1}\over{u_1\cdots u_r}}\left(\sum_{i=0}^r
(-1)^i\Bigl({{r-1}\atop{i}}\Bigr)u_{i+1}\right).\eqno(6.1)$$
This mould is actually quite similar to a power series often studied in
classical situations.  Indeed, if we define $\dar$ to be the mould operator
defined by
$$\dar\cdot A(u_1,\ldots,u_r)=u_1\cdots u_r\ A(u_1,\ldots,u_r),$$
then $\dar\cdot dupal$ is a polynomial-valued mould, so it is the image of
a power series under $ma$; explicitly
$$\dar\cdot dupal=ma\Bigl(x-{{{\rm ad}(-y)}\over{{\rm exp}({\rm ad}(-y))-1}}(x)\Bigr).$$

Ecalle gave several equivalent definitions of the key mould $pal$, but
the most recent one (see [E2]) appears to be the simplest and most convenient.
If we define $\dur$ to be the mould operator defined by
$$\dur\cdot A(u_1,\ldots,u_r)=(u_1+\cdots+u_r)\,A(u_1,\ldots,u_r),$$
then the mould $pal$ is defined recursively by 
$$\dur\cdot pal=\mu(pal,dupal).\eqno(6.2)$$
Calculating the first few terms of $pal$ explicitly, we find that
$$\cases{pal(\emptyset)=1\cr
pal(u_1)={{1}\over{2u_1}}\cr
pal(u_1,u_2)={{u_1+2u_2}\over{12u_1u_2(u_1+v_2)}}\cr
pal(u_1,u_2,u_3)={{-1}\over{24u_1u_3(u_1+u_2)}}.}$$
\vskip .3cm
Let $pil=\swap(pal)$. 
The most important result concerning $pal$, necessary for the proof of
Ecalle's fundamental identity below, is the following.
 \vskip .2cm
\noindent {\bf Theorem 6.1.} {\it The moulds $pal$ and $pil$ are symmetral.}
\vskip .2cm
In [E1,\S 4.2], the mould $pil$ (called $\frak{ess}$) is given an 
independent definition which makes it easy to prove that it is symmetral.
Similarly, it is not too hard to prove that $pal$ is symmetral
using the definition (6.2). The real difficulty is to prove that
$pil$ (as defined in [E1]) is the swap of $pal$ (as defined in (6.2)).
Ecalle sketched beautiful proofs of these two facts in [E2], 
and the details are fully written out in [S,\S\S 4.2,4.3].
\vskip .2cm
Before proceeding to the fundamental identity, we need a useful result 
in which a very simple $v$-mould is used to 
give what amounts to an equivalent definition of 
alternility.\footnote{$^*$}{This is just one example of a general identity valid for {\it flexion units}, see [E1, p. 64] where Ecalle explains the notion of
alternality twisted by a flexion unit and asserts that alternility is merely
alternality twisted by the flexion unit $1/v_1$.} 
\vskip .3cm
\noindent {\bf Proposition 6.2.} {\it Let $pic$ be the $v$-mould defined by
${\rm pic}(v_1,\ldots,v_r)=1/v_1\cdots v_r$.  Then for any alternal mould
$A\in\overline{\ARI}$, the mould $\ganit(pic)\cdot A$ is alternil.}
\vskip .3cm
\noindent Proof. The proof is deferred to Appendix C.\hfill{$\diamondsuit$}
\vskip .3cm 
We now come to Ecalle's fundamental identity.
\vskip .3cm
\noindent {\bf Ecalle's fundamental identity:} For any 
push-invariant mould $A$, we have 
$$\swap\bigl(\adari(pal)\cdot A\bigr)=\ganit(pic)\cdot \bigl(\adari(pil)\cdot \swap(A) \bigr).\eqno(6.3)$$
The proof of this fundamental identity actually follows as a consequence
of (3.2) and a more general fundamental identity, similar but taking place
in the group $\GARI$ and valid for all moulds. It is given in full detail 
in [S, Thm.~4.5.2].
\vskip .5cm
\noindent {\bf \S 7. The main theorem}
\vskip .3cm
In this section we give Ecalle's main theorem on dimorphy, which shows
how the mould $pal$ transforms moulds with the double symmetry
${\underline{al}*\underline{al}}$ to moulds that are 
${\underline{al}*\underline{il}}$. We then show how Racinet's
theorem follows directly from this.  We first need a useful lemma.
\vskip .3cm
\noindent {\bf Lemma 7.1.} {\it If $C$ is a constant-valued mould, then
$$\ganit(pic)\cdot \adari(pil)\cdot C=C.\eqno(7.1)$$}
\noindent Proof. [B, Corollary 4.43] 
We apply the fundamental identity (6.3) in the case where
$A=\swap(A)=C$ is a constant-valued mould, obtaining
$$\swap\bigl(\adari(pal)\cdot C\bigr)=\ganit(pic)\cdot \adari(pil)\cdot C.$$
So it is enough to show that the left-hand side of this is equal to $C$,
i.e.~that $\adari(pal)\cdot C=C$, since a constant mould is equal to its
own swap.  As we saw just before (3.3), the definitions (2.6)-(2.9) 
imply that $\ari(A,C)=0$ for all $A\in\ARI$. 
Now, by (5.4) we see that $\adari(pal)\cdot C$ is a linear 
combination of iterated $\ari$-brackets of $\logari(pal)$ with $C$, but
since $pal\in\GARI$, $\logari(pal)\in \ARI$, so 
$\ari(\logari(pal),C)=0$, i.e. all the bracketed terms in (5.4) are 0. 
Thus $\adari(pal)\cdot C=C$. This concludes the proof.\hfill{$\diamondsuit$}
\vskip .3cm
We can now state the main theorem on moulds.
\vskip .3cm
\noindent {\bf Theorem 7.2.} {\it The action of the operator $\adari(pal)$
on the Lie subalgebra $\ARI_{\underline{al}*\underline{al}}\subset \ARI$ yields a
Lie isomorphism of subspaces
$$\adari(pal):\ARI_{\underline{al}*\underline{al}} 
\buildrel\sim\over\rightarrow \ARI_{\underline{al}*\underline{il}}.\eqno(7.2)$$
Thus in particular $\ARI_{\underline{al}*\underline{il}}$ forms a Lie algebra
under the $\ari$-bracket.}
\vskip .2cm
\noindent Proof. 
The proof we give appears not to have been published anywhere by Ecalle, but
we learned its outline from him through a personal communication to the
second author, for which we are grateful.

Note first that $\adari(pal)$ preserves the depth 1 component of moulds
in $\ARI$, so if $A$ is even in depth 1 then so is $\adari(pal)\cdot A$.
We first consider the case where $A\in \ARI_{\underline{al}/\underline{al}}$,
i.e. $swap(A)$ is alternal without addition of a constant correction.
By (5.8), the mould $\adari(pal)\cdot A$ is alternal,
since $pal$ is symmetral by Theorem 6.1. By Proposition 3.4, $A$
is push-invariant, so Ecalle's fundamental identity (6.3)
holds.  Since $A\in \ARI_{\underline{al}/\underline{al}}$, $\swap(A)$ is 
alternal, and by Theorem 6.1, $pil$ is alternal; thus by (5.8),
$\adari(pil)\cdot \swap(A)$ is alternal. Then by Proposition 6.2, 
$\ganit(pic)\cdot \adari(pil)\cdot \swap(A)$ is alternil, and finally by
(6.3), $\swap\bigl(\adari(pal)\cdot A\bigr)$ is alternil, which proves that
$\adari(pal)\cdot A\in \ARI_{\underline{al}/\underline{il}}$ as desired.

We now consider the general case where $A\in
\ARI_{\underline{al}*\underline{al}}$.  Let $C$ be the 
constant-valued mould such that $\swap(A)+C$ is alternal. 
As above, we have that $\adari(pal)\cdot A$ is alternal, so to conclude the 
proof of the theorem it remains only to show 
that its swap is alternil up to addition of a constant mould, and
we will show that this constant mould is exactly $C$. 
As before, since $\swap(A)+C\in\overline{\ARI}$ is alternal, the mould
$$\adari(pil)\cdot \bigl(\swap(A)+C\bigr)=\adari(pil)\cdot \swap(A)+
\adari(pil)\cdot C$$ 
is also alternal.  Thus by Proposition 6.2, applying $\ganit(pic)$ to it yields 
the alternil mould
$$\ganit(pic)\cdot \adari(pil)\cdot \swap(A) + \ganit(pic)\cdot 
\adari(pil)\cdot C.$$
By Lemma 7.1, this is equal to
$$\ganit(pic)\cdot \adari(pil)\cdot \swap(A) + C,\eqno(7.3)$$
which is thus alternil. Now, since $A$ is push-invariant by Proposition 3.4,
we can apply (6.3) and find that (7.3) is equal to
$$\swap\bigl(\adari(pal)\cdot A\bigr)+C,$$
which is thus also alternil.  Therefore $\swap\bigl(\adari(pal)\cdot A\bigr)$
is alternil up to a constant, which precisely means that $\adari(pal)\cdot A
\in \ARI_{\underline{al}*\underline{il}}$ as claimed. Since $\adari(pal)$ is 
invertible (with inverse $\adari\bigl(\invgari(pal)\bigr)$) and by the analogous arguments this 
inverse takes $\ARI_{\underline{al}*\underline{il}}$ to $\ARI_{\underline{al}*
\underline{al}}$, this proves that (7.2) is an isomorphism.
\hfill{$\diamondsuit$}
\vskip .3cm
\noindent {\bf Corollary 7.3.} {\it $\ARI^{pol}_{\underline{al}*\underline{il}}$
forms a Lie algebra under the {\rm ari}-bracket.}
\vskip .2cm
\noindent Proof. By Proposition 3.2, $\ARI^{pol}$ is a Lie algebra under
the $\ari$-bracket, so since $\ARI_{\underline{al}*\underline{il}}$ is as well 
by Theorem 7.2, their intersection also forms a Lie algebra.\hfill{$\diamondsuit$}
\vskip .3cm
In view of (4.5) and (4.6), this corollary is equivalent to Racinet's
theorem that $\ds$ is a Lie algebra under the Poisson bracket.
\vfill\eject
\noindent {\bf Appendix A.}
\vskip .4cm
\noindent {\bf Proof of Theorem 3.1.} We cut it into two separate results
as explained in the main text.
\vskip .3cm
\noindent {\bf Proposition A.1.} {\it If $A$, $B$ are alternal moulds then
$C=\lu(A,B)$ is alternal.}
\vskip .2cm
\noindent Proof. We have
$$C(\w)=\lu(A,B)(\w)=\sum_{\w=\a\b} \bigl(A(\a)B(\b)-B(\a)A(\b)\bigr),$$
so we need to show that the following sum vanishes:
$$\eqalign{\sum_{\w\in\sh(\u,\v)} C(\w)&=\sum_{\w\in \sh(\u,\v)}\lu(A,B)(\w)\cr
&=\sum_{\w\in\sh(\u,\v)}\sum_{\w=\a\b} \bigl(A(\a)B(\b)-B(\a)A(\b)\bigr).}
\eqno(A.0)$$
This sum breaks into three pieces: the terms where $\a$
contain letters from both $\u$ and $\v$, the case where $\a$ contains
only letters from $\u$ or from $\v$ but $\b$ contains letters from both, 
and finally the cases $\a=\u$, $\b=\v$ and $\a=\v,\b=\u$.

The first type of terms add up to zero because we can break up the sum
into smaller sums where $\a$ lies in the shuffle of the first $i$ letters
of $\u$ and $j$ letters of $\b$, and these terms already sum to zero
since $A$ and $B$ are alternal.

The second type of term adds up to zero for the same reason, because even
though $\a$ may contain only letters from one of $\u$ and $\v$, $\b$
must contain letters from both and therefore the same reasoning holds.

The third type of term yields $A(\u)B(\v)-B(\u)A(\v)$ when $\a=\u$,$\b=\v$
and $A(\v)B(\u)-B(\v)A(\u)$ when $\a=\v$, $\b=\u$, which cancel out. Thus
the sum $(A.0)$ adds up to zero.\hfill{$\diamondsuit$}
\vskip .5cm
\noindent {\bf Proposition A.2.} {\it If $A$ and $B$ are alternal moulds in $\ARI$, 
then $C=\arit(B)\cdot A$ is alternal. }
\vskip .2cm
\noindent Proof.
By definition, $C$ is alternal if $$\sum_{\w=\sh({\bf x},{\bf y})}C(\w)=0,$$
for all pairs of non-trivial words ${\bf x},{\bf y}$. 

Pick an arbitrary pair of non-trivial words ${\bf x}, {\bf y}$, of appropriate length (that is, so that their lengths add up to the length of $A$ plus the length of $B$). We will be shuffling ${\bf x}$ and ${\bf y}$ together, and the resulting word is then broken up into three parts (all possible ones) in order to compute the flexions. Thus, if we break up $\w={\bf abc}$, ${\bf a}$ must be a shuffle of some parts at the beginning of each word ${\bf x,y}$, ${\bf b}$ must come from shuffling their middles, and ${\bf c}$ can only come from shuffling the last parts. Then we can rewrite this computation as follows:

$$\eqalign{\sum_{\w=\sh({\bf x},{\bf y})}\arit(B)\cdot A(\w)&= \sum_{\w=\sh({\bf x},{\bf y})}\left(\sum_{{{\w={\bf abc}}\atop{{\bf c}\neq\emptyset}}}A({\bf a\lceil c})B({\bf b})-\sum_{{{\w={\bf abc}\atop{\bf a}\neq\emptyset}}}A({\bf a\rceil c})B({\bf b})\right)\cr
&=\sum_{{{{\bf x=x_1x_2x_3}\atop{\bf y=y_1y_2y_3},{\bf x_3y_3\neq\emptyset}}}} \sum_{{{{\bf a}=\sh({\bf x_1},{\bf y_1})\atop{\bf b}=\sh({\bf x_2},{\bf y_2}),{\bf c}=\sh({\bf x_3},{\bf y_3})}}}A({\bf a\lceil c})B({\bf b})\cr
&-\sum_{{{{\bf x=x_1x_2x_3}\atop{\bf y=y_1y_2y_3},{\bf x_1y_1\neq\emptyset}}}} \sum_{{{{\bf a}=\sh({\bf x_1},{\bf y_1})\atop{\bf b}=\sh({\bf x_2},{\bf y_2}),{\bf c}=\sh({\bf x_3},{\bf y_3})}}}A({\bf a\rceil c})B({\bf b}).}$$

Now for a fixed splitting of each ${\bf x}$ and ${\bf y}$ into three parts, we have the following possibilities. 

\vskip .2cm
\noindent {\bf Case I.} Both ${\bf x_2}={\bf y_2}=\emptyset$. Then $B(\emptyset)=0$ so we are done. 

\vskip .2cm
\noindent {\bf Case II.} Both ${\bf x_2}$ and ${\bf y_2}$ are nonempty. The trick here is that because of the flexion operations, no matter how ${\bf b}=\sh({\bf x_2, y_2})$ is shuffled, the part being added together with the last letter in ${\bf a}$ and the first letter in ${\bf c}$ remains the same. Thus, if we further fix a particular ${\bf a}$ and ${\bf c}$, we get that 
$$\sum_{{\bf b}=\sh({\bf x_2,y_2}) }A({\bf a\lceil c})B({\bf b})=A({\bf a\lceil c})\sum_{{\bf b}=\sh({\bf x_2,y_2})}B({\bf b})=0 $$ and $$\sum_{{\bf b}=\sh({\bf x_2,y_2}) }A({\bf a\rceil c})B({\bf b})=A({\bf a\rceil c})\sum_{{\bf b}=\sh({\bf x_2,y_2})}B({\bf b})=0, $$ by alternality of $B$. And thus,

$$\sum_{{{{\bf a}=\sh({\bf x_1},{\bf y_1})\atop{\bf c}=\sh({\bf x_3},{\bf y_3})}}}\sum_{{\bf b}=\sh({\bf x_2,y_2}) }A({\bf a\lceil c})B({\bf b})=0$$ and $$\sum_{
{{{\bf a}=\sh({\bf x_1},{\bf y_1})}\atop{{\bf c}=\sh({\bf x_3},{\bf y_3})}}}\sum_{{\bf b}=\sh({\bf x_2,y_2}) }A({\bf a\rceil c})B({\bf b})=0.$$

  \noindent {\bf Case III.} Either ${\bf x_2}=\emptyset$ or ${\bf y_2}=\emptyset$, but not both. Without loss of generality, assume ${\bf x_2}=\emptyset$. Then we have:
$$\sum_{{{{\bf a}=\sh({\bf x_1},{\bf y_1})}\atop{{\bf b}={\bf y_2},{\bf c}=\sh({\bf x_3},{\bf y_3})}}}A({\bf a\lceil c})B({\bf b}) =
B({\bf y_2})\sum_{{{{\bf a}=\sh({\bf x_1},{\bf y_1})}\atop{{\bf c}=\sh({\bf x_3},{\bf y_3})}}}A({\bf a\lceil c})
$$

And similarly

$$\sum_{{{{\bf a}=\sh({\bf x_1},{\bf y_1})}\atop{{\bf b}={\bf y_2},{\bf c}=\sh({\bf x_3},{\bf y_3})}}}A({\bf a\rceil c})B({\bf b}) =B({\bf y_2})\sum_{{{{\bf a}=\sh({\bf x_1},{\bf y_1})}\atop{{\bf c}=\sh({\bf x_3},{\bf y_3})}}}A({\bf a\rceil c})
$$

Recall that by definition $$\sh({\bf x_1, y_1})=\sh({\bf x_1',y_1})(\hbox{last letter in ${\bf x_1}$})+\sh({\bf x_1, y_1'})(\hbox{last letter in ${\bf y_1}$})$$ and $$\sh({\bf x_3, y_3})=(\hbox{first letter in ${\bf x_3}$})\sh({\bf x_3',y_3})+(\hbox{first letter in ${\bf y_3}$})\sh({\bf x_3, y_3'}).$$ 

Thus, 
$${\bf a\lceil c}= \sh({\bf x_1,y_1})(\hbox{sum of letters in ${\bf y_2}$ plus first letter in ${\bf x_3}$})\sh({\bf x_3',y_3})
\eqno(A.1)$$ or 
$${\bf a\lceil c}= \sh({\bf x_1,y_1})(\hbox{sum of letters in ${\bf y_2}$ plus first letter in ${\bf y_3}$})\sh({\bf x_3,y_3'})\eqno{(A.2)}$$
and 
$${\bf a\rceil c}=\sh({\bf x_1',y_1})(\hbox{sum of letters in ${\bf y_2}$ plus last letter in ${\bf x_1}$})\sh({\bf x_3, y_3})\eqno(A.3)$$
or 
$${\bf a\rceil c}=\sh({\bf x_1,y_1'})(\hbox{sum of letters in ${\bf y_2}$ plus last letter in ${\bf y_1}$})\sh({\bf x_3, y_3}).\eqno(A.4)$$

Recall that, since ${\bf x_2}$ is assumed to be empty, then for a given ${\bf x_1, x_3}$, we can let ${\bf \overline{x_1},\overline{ x_3}}$ be so that ${\bf \overline{x_1}}$ is ${\bf x_1}$ with an additional letter given by the first letter of ${\bf x_3}$ and ${\bf \overline{x_3}}$ is defined in the logical way. That means that equations $(A.1)$ and $(A.3)$ are exactly the same. Thus, we get direct cancellation for all possible choices of ${\bf x_1, x_3}$ (this is compatible with the restrictions on nonemptiness given by the definition). 

We cannot do the same for $(A.2)$ and $(A.4)$, since ${\bf y_2}$ is assumed to be nonempty. For these, notice that if we keep ${\bf y}$ fixed and sum over all possible partitions of ${\bf x=x_1x_2x_3}$ where ${\bf x_2}=\emptyset$, and ${\bf x_3}\neq\emptyset$ we get that each
$${\bf a\lceil c}= \sh({\bf x_1,y_1})(\hbox{sum of letters in ${\bf y_2}$ plus first letter in ${\bf y_3}$})\sh({\bf x_3,y_3'})$$ 
could be seen as a term in the shuffle $\sh({\bf x,y_1\lceil y_3})$. To see this, suppose that $${\bf x}=u_1\cdots u_k|u_{k+1}\cdots u_l={\bf x_1|x_3}$$ and that $${\bf y}=u_{l+1}\cdots u_{l+i}|u_{l+i+1}\cdots u_{l+j}|u_{l+j+1}\cdots u_n={\bf y_1|y_2|y_3}.$$ 

Then $${\bf a\lceil c}=\sh((u_1\cdots u_k),(u_{l+1}\cdots u_{l+i}))(u_{l+i+1}+\cdots+u_{l+j}+u_{l+j+1})\sh((u_{k+1}\cdots u_l),(u_{l+j+2}\cdots u_n)).$$

And so if we allow the $k$ to shift from 1 to $l$, this is essentially the shuffling of the words $u_1\cdots u_l={\bf x}$ and $u_{l+1}\cdots u_{l+i}(u_{l+i+1}+\cdots+u_{l+j}+u_{l+j+1})u_{l+j+2}\cdots u_n={\bf y_1\lceil y_3}$. Thus we have

$$\sum_{{{{\bf x=x_1x_3}}\atop{{\bf x_3}\neq\emptyset}}} \sum_{{{{\bf a}=\sh({\bf x_1},{\bf y_1})}\atop{{\bf b}={\bf y_2},{\bf c}=y_{\hbox{first}}\sh({\bf x_3},{\bf y_3'})}}}A({\bf a\lceil c})=\sum_{\w=\sh({\bf x, y_1\lceil y_3})}A(\w)=0$$ by alternality of $A$. 

A similar argument holds for the terms corresponding to the other flexion (the terms corresponding to equation $(A.4))$. 

Putting all of these cases together, we see that indeed, $C$ is alternal.\hfill{$\diamondsuit$} 
\vskip .2cm

\noindent {\bf Proposition A.3.} {\it If $A$ and $B$ are alternal moulds in $\overline{\ARI}$, 
then $C=\arit(B)\cdot A$ is alternal. }
\vskip .2cm
\noindent Proof.  As with the proof for $\ARI_{al}$, we have to show that $$\sum_{\w=\sh({\bf x},{\bf y})}C(\w)=0,$$
for all pairs of non-trivial words ${\bf x},{\bf y}$. Again, this can be rewritten as follows:

$$\eqalign{\sum_{\w=\sh({\bf x},{\bf y})}\arit(B)\cdot A(\w)&= \sum_{\w=\sh({\bf x},{\bf y})}\left(\sum_{\w={\bf abc}\atop{\bf c}\neq\emptyset}A({\bf a c})B({\bf b}\rfloor)-\sum_{\w={\bf abc}\atop{\bf a}\neq\emptyset}A({\bf a c})B({\bf \lfloor b})\right)\cr
&=\sum_{{\bf x=x_1x_2x_3}\atop{\bf y=y_1y_2y_3},{\bf x_3y_3\neq\emptyset}} \sum_{{\bf a}=\sh({\bf x_1},{\bf y_1})\atop{\bf b}=\sh({\bf x_2},{\bf y_2}),{\bf c}=\sh({\bf x_3},{\bf y_3})}A({\bf a c})B({\bf b\rfloor})\cr
&-\sum_{{\bf x=x_1x_2x_3}\atop{\bf y=y_1y_2y_3},{\bf x_1y_1\neq\emptyset}} \sum_{{\bf a}=\sh({\bf x_1},{\bf y_1})\atop{\bf b}=\sh({\bf x_2},{\bf y_2}),{\bf c}=\sh({\bf x_3},{\bf y_3})}A({\bf a c})B({\bf \lfloor b})}$$

Again, for a fixed splitting of each ${\bf x}$ and ${\bf y}$ into three parts, we have the following possibilities. 
\vskip .2cm
\noindent {\bf Case I.} Both ${\bf x_2}={\bf y_2}=\emptyset$. Then $B(\emptyset)=0$ so we are done. 
\vskip .2cm
\noindent {\bf Case II.} Both ${\bf x_2}$ and ${\bf y_2}$ are nonempty. 

Here, no matter how ${\bf b}=\sh({\bf x_2, y_2})$ is shuffled, the part being subtracted from ${\bf b}$, which is either the last letter in ${\bf a}$ or the first letter in ${\bf c}$, remains the same if we fix a particular ${\bf a}$ and ${\bf c}$. Thus, we get that  
$${\bf b\rfloor}_i=\sh({\bf x_2, y_2})_i-\hbox{first letter in ${\bf c}$}=\sh(({\bf x_{2_k}}-\hbox{first letter in ${\bf c}$}),({\bf y_{2_k}}-\hbox{first letter in ${\bf c}$}))_i$$ and $${\bf \lfloor b}_i=\sh({\bf x_2, y_2})_i-\hbox{last letter in ${\bf a}$}=\sh(({\bf x_{2_k}}-\hbox{last letter in ${\bf a}$}),({\bf y_{2_k}}-\hbox{last letter in ${\bf a}$}))_i.$$

Thus,

$$\sum_{{\bf b}=\sh({\bf x_2,y_2}) }A({\bf a c})B({\bf b\rfloor})=A({\bf a c})\sum_{{\bf b}=\sh({\bf x_2,y_2})}B({\bf b\rfloor})=0 $$ and $$\sum_{{\bf b}=\sh({\bf x_2,y_2}) }A({\bf a c})B({\bf \lfloor b})=A({\bf a c})\sum_{{\bf b}=\sh({\bf x_2,y_2})}B({\bf \lfloor b})=0, $$ by alternality of $B$. And thus,

$$\sum_{{\bf a}=\sh({\bf x_1},{\bf y_1})\atop{\bf c}=\sh({\bf x_3},{\bf y_3})}\sum_{{\bf b}=\sh({\bf x_2,y_2}) }A({\bf a c})B({\bf b\rfloor})=0$$ and $$\sum_{{\bf a}=\sh({\bf x_1},{\bf y_1})\atop{\bf c}=\sh({\bf x_3},{\bf y_3})}\sum_{{\bf b}=\sh({\bf x_2,y_2}) }A({\bf a c})B({\bf \lfloor b})=0.$$

  \noindent {\bf Case III.} Either ${\bf x_2}=\emptyset$ or ${\bf y_2}=\emptyset$, but not both. Without loss of generality, assume ${\bf x_2}=\emptyset$. 
  
  Recall, again, that by definition $$\sh({\bf x_1, y_1})=\sh({\bf x_1',y_1})(\hbox{last letter in ${\bf x_1}$})+\sh({\bf x_1, y_1'})(\hbox{last letter in ${\bf y_1}$})$$ and $$\sh({\bf x_3, y_3})=(\hbox{first letter in ${\bf x_3}$})\sh({\bf x_3',y_3})+(\hbox{first letter in ${\bf y_3}$})\sh({\bf x_3, y_3'}).$$ 
  
  Since ${\bf x_2}=\emptyset$, we can see that 
  
  $${\bf b\rfloor}_i= {\bf y_{2_i}}-\hbox{first letter in ${\bf c}$}$$ and $${\bf \lfloor b}_i={\bf y_{2_i}}-\hbox{last letter in ${\bf a}$}.$$
  
 For a given ${\bf x_1, x_3}$, we can let ${\bf \overline{x_1},\overline{ x_3}}$ be so that ${\bf \overline{x_1}}$ is ${\bf x_1}$ with an additional letter given by the first letter of ${\bf x_3}$ and ${\bf \overline{x_3}}$ is defined in the logical way. That means that $$A(\sh({\bf \overline{x_1}',y_1})(\hbox{last letter in ${\bf \overline{x_1}}$})\sh({\bf \overline{x_3}, y_3}))B({\bf \lfloor b})$$ and $$A(\sh({\bf {x_1},y_1})(\hbox{first letter in ${\bf x_3}$})\sh({\bf x_3', y_3}))B({\bf \b\rfloor})$$ are identical (for each fixed shuffling).

 Thus, we get direct cancellation for all possible choices of ${\bf x_1, x_3}$ (this is compatible with the restrictions on nonemptiness given by the definition). 
 
 The only terms that have not cancelled out are the ones coming from the second term in the shuffle equations above. Now, suppose that $${\bf x}=v_1\cdots v_k|v_{k+1}\cdots v_l={\bf x_1|x_3}$$ and that $${\bf y}=v_{l+1}\cdots v_{l+i}|v_{l+i+1}\cdots v_{l+j}|v_{l+j+1}\cdots v_n={\bf y_1|y_2|y_3}, $$ and fix this splitting of ${\bf y}$. Then $${\bf ac}=\sh({ v_1\cdots v_k, v_{l+1}\cdots v_{l+i}})v_{l+j+1}\sh({ v_{k+1}\cdots v_l, v_{l+j+2}\cdots v_n}). $$  
 And so if we allow the $k$ to shift from 1 to $l$, this is essentially the shuffling of the words $v_1\cdots v_l={\bf x}$ and $v_{l+1}\cdots v_{l+i},v_{l+j+1},v_{l+j+2}\cdots v_n={\bf y_1 y_3}$. Notice that this shuffling fixes ${\bf b\rfloor}$, since $${\bf b\rfloor}=(v_{l+i+1}-v_{l+j+1},\dots,v_{l+j}-v_{l+j+1}).$$
 Thus we have

$$\sum_{{\bf x=x_1x_3}\atop{\bf x_3}\neq\emptyset} \sum_{{\bf a}=\sh({\bf x_1},{\bf y_1}),{\bf b}={\bf y_2}\atop{\bf c}=y_{\hbox{first}}\sh({\bf x_3},{\bf y_3'})}A({\bf a c})B({\bf b\rfloor})=B({\bf b\rfloor})\sum_{w=\sh({\bf x, y_1 y_3})}A(w)=0$$ by alternality of $A$. 

A similar argument holds for the terms corresponding to the other flexion.  Combining all the cases, we see that indeed, $C$ is alternal. \hfill{$\diamondsuit$}

\vfill\eject
\noindent {\bf Appendix B.}
\vskip .4cm
\noindent {\bf Proof of Proposition 3.4.}  
By additivity, we may assume that $A$ is concentrated in a fixed
depth $d$, meaning that  $A(u_1,\ldots,u_r)=0$ for all $r\ne d$.  
We use the following two lemmas.
\vskip .3cm
\noindent {\bf Lemma B.1.} {\it If $A\in \ARI_{al}$, then
$$A(u_1,\ldots,u_r)=(-1)^{r-1}A(u_r,\ldots,u_1);$$
in other words, $A$ is $\mantar$-invariant. Similarly, if $A\in\overline{\ARI}_{al}$
then again $A$ is $\mantar$-invariant}
\vskip .3cm
\noindent Proof. We give the argument for $\ARI$; the result in $\overline{\ARI}$
comes from the identical computation with $u_i$ replaced by $v_i$.
We first show that the sum of shuffle relations
$$\sh\bigl((1),(2,\ldots,r)\bigr)-\sh\bigl((2,1),(3,\ldots,r)\bigr)
+\sh\bigl((3,2,1),(4,\ldots,r)\bigr)+\cdots$$
$$\ \ \ \ \ +(-1)^{r-1}\sh\bigl((r-1,\ldots,2,1),(r)\bigr)=(1,\ldots,r)+(-1)^{r-1}
(r,\ldots,1).$$
Indeed, using the recursive formula for shuffle, we can write the above sum with two terms for each shuffle, as
$$\eqalign{(1,\ldots,r)&+2\cdot \sh\bigl((1),(3,\ldots,r)\bigr)\cr
&-2\cdot \sh\bigl((1),(3,\ldots,r)\bigr)-3\cdot \sh\bigl((2,1),(4,\ldots,r)\bigr)\cr
&+3\cdot \sh\bigl((2,1),(4,\ldots,r)\bigr)+4\cdot \sh\bigl((3,2,1),(5,\ldots,r)\bigr)\cr
&+\cdots+(-1)^{r-2}(r-1)\cdot \sh\bigl((r-2,\ldots,1),(r)\bigr)\cr
&+(-1)^{r-1}(r-1)\cdot \sh\bigl((r-2,\ldots,1),(r)\bigr)+(-1)^{r-1}(r,r-1,\ldots,1)\cr
&=(1,\ldots,r)+(-1)^{r-1}(r,\ldots,1).}$$
Using this, we conclude that if $A$ satisfies the shuffle relations, then 
$$A(u_1,\ldots,u_r)+(-1)^{r-1}A(u_r,\ldots,u_1),$$
which is the desired result.\hfill{$\diamondsuit$}
\vskip .3cm
\noindent {\bf Lemma B.2.} {\it If $A\in \ARI_{\underline{al}*\underline{al}}$, 
then $A$ is $\neg\circ \push$-invariant.}
\vskip .3cm
\noindent Proof.  We first consider the case where 
$A\in \ARI_{\underline{al}/\underline{al}}$. Using the easily verified identity
$$\neg\circ \push=\mantar\circ \swap\circ \mantar\circ \swap,\eqno(B.1)$$
and the fact that by Lemma B.1, if $A\in\ARI_{\underline{al}/\underline{al}}$, then both
$A$ and $\swap(A)$ are $\mantar$-invariant, we have
$$\eqalign{\neg\circ \push(A)(u_1,\ldots,u_r)&=
\mantar\circ \swap\circ \mantar\circ \swap(A)(u_1,\ldots,u_r)\cr
&=(-1)^{r-1}\mantar\circ \swap\circ \swap(A)(u_1,\ldots,u_r)\cr
&=(-1)^{r-1}\mantar(A)(u_1,\ldots,u_r)\cr
&=A(u_1,\ldots,u_r),}\eqno(B.2)$$
so $A$ is $\neg\circ \push$-invariant.

Now suppose that $A\in\ARI_{\underline{al}*\underline{al}}$, so $A$ is
alternal and $\swap(A)+A_0$ is alternal for some constant mould $A_0$.
By additivity, we may assume that $A$ is concentrated in depth $r$.  First
suppose that $r$ is odd.  Then $\mantar(A_0)(v_1,\ldots,v_r)=(-1)^{r-1}
A_0(v_r,\ldots,v_1)$, so since $A_0$ is a constant mould, it is 
mantar-invariant.  But $\swap(A)+A_0$ is alternal, so it is also 
mantar-invariant by Lemma B.1; thus $\swap(A)$ is mantar-invariant, and
the identity $\neg\circ\push=\mantar\circ\swap\circ\mantar\circ\swap$
shows that $A$ is $\neg\circ\push$-invariant as in $(B.2)$.

Finally, we assume that $A$ is concentrated in even depth $r$. Here
we have $\mantar(A_0)=-A_0$, so we cannot use the argument above;
indeed $\swap(A)+A_0$ is mantar-invariant, but 
$$\mantar(\swap(A))=\swap(A)+2A_0.\eqno(B.3)$$  
Instead, we note that if $A$ is alternal then so is $\neg(A)=A$.  Thus we can
write $A$ as a sum of an even and an odd function of the $u_i$ via the formula
$$A={{1}\over{2}}(A+\neg(A))+{{1}\over{2}}(A-\neg(A)).\eqno(B.4)$$
So it is enough to prove the desired result for all moulds concentrated in
even depth $r$ such that either $\neg(A)=A$ (even functions) or
$\neg(A)=-A$ (odd functions).  First suppose that $A$ is
even. Then since $\neg$ commutes with $\push$ and $\push$ is of odd order
$r+1$ and $\neg$ is of order 2, we have 
$$(\neg\circ\push)^{r+1}(A)=\neg(A)=A.\eqno(B.5)$$  
However, we also have 
$$\eqalign{\neg\circ\push(A)&=\mantar\circ\swap\circ\mantar\circ\swap(A)\cr
&=\mantar\circ\swap\bigl(\swap(A)+2A_0\bigr)\ \ {\rm by \ } (B.3)\cr
&=\mantar\bigl(A+2A_0\bigr)\cr
&=A-2A_0.}$$
Thus $(\neg\circ\push)^{r+1}(A)=A-2(r+1)A_0$, and this is equal to $A$
by $(B.5)$, so $A_0=0$; thus in fact $A\in\ARI_{\underline{al}/\underline{al}}$ 
and that case is already proven. 

Finally, if $A$ is odd, i.e. $\neg(A)=-A$, the same argument as above
gives $A-2(r+1)A_0=-A$, so $A=(r+1)A_0$, so $A$ is a constant-valued
mould concentrated in depth $r$, but this contradicts the assumption
that $A$ is alternal since constant moulds are not alternal, unless 
$A=A_0=0$.  Note that this argument shows that all moulds in 
$\ARI_{\underline{al}*\underline{al}}$ that are not in 
$\ARI_{\underline{al}/\underline{al}}$ must be concentrated in odd depths.
\hfill{$\diamondsuit$}
\vskip .4cm
We can now complete the proof of Proposition 3.4\footnote{$^*$}{Ecalle
states this result in [E1,\S 2.4] and there is also a proof in
the preprint [E2,\S 12], but we were not able to follow the argument, 
so we have provided this alternative proof.}.
Because $A=\neg\circ \push(A)$, we have $\neg(A)=\push(A)$, so in fact
we only need to show that $\neg(A)=A$.  As before, we may assume that
$A$ is concentrated in depth $r$. If $r=1$, then $A$ is an even function
by assumption.  If $r$ is even, then as before we have
$A=(\neg\circ \push)^{2s+1}(A)=\neg(A).$
Finally, assume $r=2s+1$ is odd. Since we can write $A$ as a sum of an 
even and an odd part as in (B.4),
we may assume that $\neg(A)=-A$.
Then, since $A$ is alternal, using the shuffle
$\sh\bigl((u_1,\ldots,u_{2s})(u_{2s+1})\bigr)$, we have
$$\sum_{i=0}^{2s}A(u_1,\ldots,u_i,u_{2s+1},u_{i+1},\ldots,u_{2s})=0.$$
Making the variable change $u_0\leftrightarrow u_{2s+1}$ gives
$$\sum_{i=0}^{2s}A(u_1,\ldots,u_i,u_0,u_{i+1},\ldots,u_{2s})
=0.\eqno(B.6)$$
Now consider the shuffle relation $\sh((u_1)(u_2,\ldots,u_{2s+1}))$, which
gives
$$\sum_{i=1}^{2s+1} A(u_2,\ldots,u_i,u_1,u_{i+1},\ldots,u_{2s+1})=0.\eqno(B.7)$$
Set $u_0=-u_1-\cdots -u_{2s+1}$.  Since $\neg\circ \push$ acts
like the identity on $A$, we can apply it to each term of $(B.7)$ to obtain
$$\sum_{i=1}^{2s} -A(u_0,u_2,\ldots,u_i,u_1,u_{i+1},\ldots,u_{2s})
-A(u_0,u_2,\ldots,u_{2s},u_{2s+1}).$$
We apply $\neg\circ \push$ again to the final term of this sum in order to
get the $u_{2s+1}$ to disappear, obtaining
$$\sum_{i=1}^{2s} -A(u_0,u_2,\ldots,u_i,u_1,u_{i+1},\ldots,u_{2s})
+A(u_1,u_0,u_2,\ldots,u_{2s-1},u_{2s})=0.$$
\vskip .1cm
Making the variable change $u_0\leftrightarrow u_1$ in this identity yields
$$\sum_{i=1}^{2s} -A(u_1,u_2,\ldots,u_i,u_0,u_{i+1},\ldots,u_{2s})
+A(u_0,u_1,u_2,\ldots,u_{2s-1},u_{2s})=0.\eqno(B.8)$$
Finally, adding $(B.6)$ and $(B.8)$ yields 
$2A(u_0,u_1,\ldots,u_{2s})=0,$
so $A=0$.  This concludes the proof that $\neg(A)=A$ for all
$A\in \ARI_{\underline{al}*\underline{al}}$, and thus, by Lemma B.2, that
$\push(A)=A$. This concludes the proof of Proposition 3.4.\hfill{$\diamondsuit$}
\vfill\eject
\noindent {\bf Appendix C.}
\vskip .4cm
We follow Ecalle's more general construction of {\it twisted alternality}
from [E1, pp.~57-64].
Let ${\bf e}\in \overline{\ARI}$ be a {\it flexion unit}, which is a mould 
concentrated in depth 1 satisfying 
$${\bf e}(v_1)=-{\bf e}(-v_1)$$
and
$${\bf e}(v_1){\bf e}(v_2)={\bf e}(v_1-v_2){\bf e}(v_2)+{\bf e}(v_1){\bf e}(v_2-v_1).$$
Associate to ${\bf e}$ the mould ${\bf ez}\in\GVARI$ defined by 
$${\bf ez}(v_1,\ldots,v_r)={\bf e}(v_1)\cdots {\bf e}(v_r).$$
Then a mould $A\in\overline{\ARI}$ is said to be {\it ${\bf e}$-alternal} if 
$A=\ganit({\bf ez})\cdot B$ where $B\in\overline{\ARI}$ is alternal.
The conditions for ${\bf e}$-alternality can be written out using the 
explicit expression for $\ganit$, using flexions, computed by Ecalle 
[E1, (2.36)]:
$$\bigl(\ganit(B)\cdot A\bigr)({\bf w})=\sum A({\bf b}^1\cdots
{\bf b}^s)B(\lfloor{\bf c}^1)\cdots A(\lfloor{\bf c}^s),\eqno(C.1)$$
where the sum runs over the decompositions of the word 
${\bf w}=(u_1,\ldots,u_r)$ ($r\ge 1$) as 
$${\bf w}={\bf b}^1{\bf c}^1\cdots {\bf b}^s{\bf c}^s,\ \ (s\ge 1)$$
where all ${\bf b}^i$ and ${\bf c}^i$ are non-empty words except
possibly for ${\bf c}^s$.  For example in small depths, setting 
$C=\ganit(B)\cdot A$, we have 
$$\cases{C(v_1)=A(v_1)\cr
C(v_1,v_2)=A(v_1,v_2)+A(v_1)B(v_2-v_1)\cr
C(v_1,v_2,v_3)=A(v_1,v_2,v_3)+A(v_1,v_2)B(v_3-v_2)+\cr
\ \qquad\qquad\qquad A(v_1)B(v_2-v_1,v_3-v_1)+A(v_1,v_3)B(v_2-v_1).}$$
Using the expression (C.1) for $\ganit(B)\cdot A$, the ${\bf e}$-alternality 
relations can be written explicitly as follows.  Let 
$Y_1=(y_1,\ldots,y_r)$ and $Y_2=(y_{r+1},\ldots,y_{r+s})$.  Then for
each word in the stuffle set $st(Y_1,Y_2)$, we construct the associated
{\it ${\bf e}$-alternality term}, with an expression of the form 
$$\bigl(C(\ldots,v_i,\ldots)-C(\ldots,v_j)\bigr){\bf e}(v_i-v_j)$$
corresponding each contraction (cf. (2.13).  For example, taking
$Y_1=(y_i,y_j)$ and $Y_2=(y_k,y_l)$, the stuffle set $st(Y_1,Y_2)$
is given in (2.12), and the corresponding 13 ${\bf e}$-alternality terms 
are, first of all the six shuffle terms
$$C(v_1,v_2,v_3,v_4),C(v_1,v_3,v_2,v_4),C(v_1,v_3,v_4,v_2),
C(v_3,v_1,v_2,v_4),$$ $$C(v_3,v_1,v_4,v_2),C(v_3,v_4,v_1,v_2)$$
(cf. (2.14)), then the six terms with a single contraction
$$\bigl(C(v_1,v_2,v_4)-C(v_1,v_3,v_4)\bigr){\bf e}(v_2-v_3),\ \ \
\bigl(C(v_1,v_2,v_4)-C(v_3,v_2,v_4)\bigr){\bf e}(v_1-v_3),$$
$$\bigl(C(v_1,v_3,v_2)-C(v_1,v_3,v_4)\bigr){\bf e}(v_2-v_4),\ \ \
\bigl(C(v_1,v_4,v_2)-C(v_3,v_4,v_2)\bigr){\bf e}(v_1-v_3),$$
$$\bigl(C(v_3,v_1,v_2)-C(v_3,v_1,v_4)\bigr){\bf e}(v_2-v_4),\ \ \
\bigl(C(v_3,v_1,v_2)-C(v_3,v_4,v_2)\bigr){\bf e}(v_1-v_4)$$
(cf. (2.15)),
and finally the single term with two contractions,
$$\bigl(C(v_1,v_2)-C(v_3,v_2)-C(v_1,v_4)+C(v_2,v_4)\bigr){\bf e}(v_1-v_3)
{\bf e}(v_2-v_4).$$
The {\it ${\bf e}$-alternality sum} $C_{r,s}$ is defined to be the sum of 
all the ${\bf e}$-alternality terms corresponding to words in the stuffle
set $st(Y_1,Y_2)$; this sum is independent of the actual sequences
$Y_1,Y_2$, depending only on their lengths $r,s$.  The mould $C$ is
said to satisfy the {\it ${\bf e}$-alternality relations} if
$C_{r,s}=0$ for all $1\le r\le s$.  Comparing with (2.14-15) we see
that the notion of alternality is nothing but the special case of
${\bf e}$-alternality for the flexion unit ${\bf e}(v_1)=1/v_1$.
The associated mould ${\bf ez}$ is thus equal to $pic$, so we
find that $\ganit(pic)\cdot A$ is alternil if $A$ is alternal.

\vskip .3cm
\vfill\eject
\noindent {\bf References}
\vskip .3cm
\noindent [B] S. Baumard, Aspects modulaires et elliptiques des relations
entre multiz\^etas, Ph.D. dissertation, Paris, France, 2014.
\vskip .3cm
\noindent [C] S. Carr, Multizeta Values: Lie algebras and periods on
$M_{0,n}$, Ph.D. Thesis, 2008.  Available at http://math.unice.fr/~brunov/GdT/Carr.pdf.
\vskip .3cm
%\noindent [Cr] J. Cresson, Calcul Moulien, {\it Ann. Fac. Sci. Toulouse Math. } {\bf XVIII} (2009), no. 2, 307-395.
%\vskip .3cm
\noindent [E1] J. Ecalle, The flexion structure and dimorphy: flexion units,
singulators, generators, and the enumeration of multizeta irreducibles,
in {\it Asymptotics in Dynamics, Geometry and PDEs; Generalized Borel Summation},
vol. II, O. Costin, F. Fauvet, F. Menous, D. Sauzin, eds., Edizioni della Normale, Pisa, 2011.
\vskip .3cm
\noindent [E2] J. Ecalle, Eupolars and their bialternality grid, preprint 2014.
\vskip .3cm
\noindent [F] H. Furusho, Double shuffle relation for associators, {\it Ann. of Math (2)} {\bf 174} (2011), no. 1, 341-360.
\vskip .3cm
\noindent [R1] G. Racinet, S\'eries g\'en\'eratrices non-commutatives de
polyz\^etas et associateurs de Drinfel'd, Ph.D. dissertation, Paris, France, 2000.
\vskip .3cm
\noindent [R2] G. Racinet, Doubles m\'elanges des polylogarithmes multiples
aux racines de l'unit\'e, {\it Publ. Math. Inst. Hautes Etudes Sci.} {\bf No. 95} (2002), 185-231.
\vskip .3cm
\noindent [S] L. Schneps, ARI, GARI, Zig and Zag: An introduction to Ecalle's
theory of multizeta values, arXiv:1507.01534, 2014.
\vskip .3cm
\noindent [Se] J-P. Serre, Lie Algebras and Lie Groups, Springer Lecture Notes
{\bf 1500}, Springer-Verlag, 1965.
\bye